\documentclass[11pt]{article}
\newtheorem{theorem}{Theorem}[section]
\newtheorem{proposition}[theorem]{Proposition}
\newtheorem{lemma}[theorem]{Lemma}

\newcommand{\bthm}{\begin{theorem}}
\newcommand{\ethm}{\end{theorem}}
\newcommand{\bprop}{\begin{proposition}}
\newcommand{\eprop}{\end{proposition}}
\newcommand{\blem}{\begin{lemma}}
\newcommand{\elem}{\end{lemma}}
\newcommand{\beq}{\begin{equation}}
\newcommand{\eeq}{\end{equation}}

\newcommand{\nl}{\newline}

\newcommand{\dist}{{\rm dist}}

\newcommand{\R}{{\bf R}}

\newcommand{\diam}{{\rm diam}}

\newcommand{\ia}{({\rm i})}
\newcommand{\ib}{({\rm ii})}
\newcommand{\ic}{({\rm iii})}
\newcommand{\hk}{\Bigl(\frac{p-1}{p}\Bigr)}

\title{Refined geometric $L^p$  Hardy inequalities}
\author{
G. Barbatis\footnote{Department of Mathematics,
 University of Ioannina, 45110 Ioannina, Greece}
 \and S. Filippas\footnote{Department of Applied Mathematics,
 University of Crete, 71409 Heraklion, Greece} \and A. Tertikas
 \footnote{Department of Mathematics,
 University of Crete, 71409 Heraklion, Greece and  \nl
Institute of  Applied and Computational Mathematics,
FORTH, 71110 Heraklion, Greece}
}
\begin{document}
\maketitle

\begin{abstract}

    For a bounded convex domain $\Omega$  in $\R^N$ we prove  refined Hardy inequalities
that involve the Hardy potential corresponding to the  distance to the boundary of $\Omega$,
 the volume of  $\Omega$, as well as
a finite number of sharp logarithmic corrections.
We also discuss the best constant of these inequalities.

\end{abstract}

\section{Introduction}
For a convex domain $\Omega\subset\R^N$ the Hardy inequality
\beq
\int_{\Omega}|\nabla u|^pdx\geq\hk^p\int_{\Omega}\frac{|u|^p}{d^p}dx,\quad d(x)=\dist(x,\partial\Omega)
\quad u\in W^{1,p}_0(\Omega)
\label{eq:hi}
\eeq
is valid, where the constant $\hk^p$ is optimal; cf [MMP], [MS].
 Brezis and Marcus  \cite{BM} have  established an improved version of (\ref{eq:hi}) when $p=2$: they showed
that for bounded and convex $\Omega$ there holds
\beq
\int_{\Omega}|\nabla u|^2dx\geq\frac{1}{4}\int_{\Omega}\frac{u^2}{d^2}dx +\frac{1}{4\diam^2(\Omega)}\int_{\Omega}u^2dx,
\qquad u\in H^1_0(\Omega).
\label{eq:bm}
\eeq
The question was asked in that paper as to whether it is possible to replace $\diam^{-2}(\Omega)$
 by $c|\Omega|^{-2/N}$, where
$|\Omega|$ denotes the volume of $\Omega$.
 A positive answer was given by  M. and T.
 Hoffmann-Ostenhof and   Laptev [HHL], who showed that
\beq
\int_{\Omega}|\nabla u|^2dx\geq\frac{1}{4}\int_{\Omega}\frac{u^2}{d^2}dx
 +k_2 \Big(\frac{a_N}{|\Omega|} \Big)^{\frac{2}{N}}    \int_{\Omega}u^2dx,
\qquad u\in H^1_0(\Omega)
\label{eq:l}
\eeq
where
where $a_N$ is the volume of the unit ball and $k_2 = N/4$.
%\[\hspace{3cm}k_2=\frac{Na_N^{2/N}}{4}\qquad\quad\mbox{($a_N=$ volume of unit ball)}\]

 In connection with this let us notice that when we take as 
$d(x)$ the distance from a point of $\Omega$, say the origin, the following improved Hardy inequality
was established by Brezis and Vazquez [BV]
\beq
\int_{\Omega}|\nabla u|^2dx\geq \Big(\frac{N-2}{2}\Big)^2 
\int_{\Omega}\frac{u^2}{|x|^2}dx + \mu_2 \Big(\frac{a_N}{|\Omega|} \Big)^{\frac{2}{N}} 
\int_{\Omega}u^2dx,
\qquad u\in H^1_0(\Omega);
\label{eq:bv}
\eeq
here  $\mu_2\simeq 5.783$ is  the first eigenvalue of the Dirichlet Laplacian for the unit disk in $\R^2$.
This constant is optimal when $\Omega$ is a ball centered at the origin, independently of
 the dimension $N \geq 2$, cf [BV],
whereas for general $\Omega$ this constant is not optimal, cf [FT, Proposition 5.1].

%\[\hspace{3cm}k_2=\frac{Na_N^{2/N}}{4}\qquad\quad\mbox{($a_N=$ volume of unit ball)}\]
An $L^p$ -version of (\ref{eq:l}) was recently obtained by Tidblom \cite{T} who showed that for
 convex $\Omega$ there holds
\beq
\int_{\Omega}|\nabla u|^pdx\geq\hk^p\int_{\Omega}\frac{|u|^p}{d^p}dx +k_p \Big(\frac{a_N}{|\Omega|} \Big)^{\frac{p}{N}}    \int_{\Omega}u^2dx,
\qquad u\in W^{1,p}_0(\Omega)
\label{eq:t}
\eeq
with
\beq
k_p=(p-1)\hk^p \frac{\sqrt{\pi}\Gamma(\frac{N+p}{2})}{\Gamma(\frac{p+1}{2})
\Gamma(\frac{N}{2})}.
\label{1.kp}
\eeq
For $p=2$ this reduces to (\ref{eq:l}); in particular $k_2 = N/4$.

In addition to (\ref{eq:l})  it was shown in [HHL, Theorem 3.4] that if 
\beq
X_1(t) = (1- \log t)^{-1},  \quad  \quad \quad   t\in (0,1),
\label{eq:x1}
\eeq
the following more refined improvement of  (\ref{eq:l}) is  true:  for any $D\geq\diam(\Omega)/2$ there holds
\begin{eqnarray}
&&\hspace{-1cm}\int_{\Omega}|\nabla u|^2dx \geq\frac{1}{4}\int_{\Omega}\frac{u^2}{d^2}dx +
 \frac{1}{4}\int_{\Omega}\frac{u^2}{d^2}X_1^2(d/D)dx \nonumber \\
&&\qquad \qquad +k_2(1-X_1(\diam(\Omega)/ (2D))^2 \left(\frac{a_N}{|\Omega|}\right)^{2/N}\int_{\Omega}u^2dx
\label{eq:mus}
\end{eqnarray}
for all $u\in H^1_0(\Omega)$. Note that if we let $D\to\infty$ in (\ref{eq:mus}) we regain (\ref{eq:l}).

In our main result we extend both (\ref{eq:t}) and (\ref{eq:mus}).
More precisely, with $X_1(t)$ as in (\ref{eq:x1})  we define recursively
\begin{equation}
 X_k(t) = X_{1}(X_{k-1}(t)), \quad k=2,3,\ldots ,t\in (0,1).
\label{eq:flash}
\end{equation}
These are iterated logarithmic functions that vanish at an increasingly low rate at $t=0$. Let us fix
$k \geq 1$ and set
\beq
a=\left\{
\begin{array}{ll}
 \mbox{$0$,} &   \mbox{ \quad if \quad  $1 < p \leq 2$,} \\
   \mbox{$ \frac{(p-2) k}{3(p-1)}>0$,}  &     \mbox{  \quad if \quad 
\quad \quad  $p>2$,}
\end{array} \right.
\label{2.a}
\eeq
and
\[\eta(t)=\sum_{i=1}^kX_1(t)\ldots X_i(t),\]
whereas for $k=0$ we set $\eta  =0$.
For $D\geq\diam(\Omega)/2$ we also set
\[\eta_D=\eta(\frac{\diam(\Omega)}{2D}).\]
Then our main result reads:

{\bf Theorem  A} {\em
Assume that $\Omega$ is convex and bounded. Let $k\geq 0$ be a fixed
integer. 
Then, there exists $D_0=  D_0(k, p ,\diam(\Omega)) \geq \diam(\Omega)/2$ such that for $D\geq D_0$ there holds
\begin{eqnarray}
&&\hspace{-1cm}\int_{\Omega}|\nabla u|^pdx \geq\hk^p\int_{\Omega}\frac{|u|^p}{d^p}dx +
 \frac{1}{2}\hk^{p-1}\sum_{i=1}^k\int_{\Omega}\frac{|u|^p}{d^p}X_1^2(d/D)\ldots X_i^2(d/D)dx \nonumber \\
&&\qquad \qquad +k_p(1-\eta_D-a\eta_D^2)^{\frac{p}{p-1}}  \Big(\frac{a_N}{|\Omega|} \Big)^{\frac{p}{N}}  \int_{\Omega}|u|^pdx, \label{eq:dos}
\end{eqnarray}
for all $u\in W^{1,p}_0(\Omega)$.
 When $p=2$ we can take as $D_0$ the unique solution  of  $\eta_{D_0}=1$.
}

Note that if we let $D\to+\infty$ in (\ref{eq:dos}) we recover (\ref{eq:t}).
Also, for $p=2$ and $k=1$ we recover (\ref{eq:mus}).
Moreover, the terms in the series are sharp: it was shown in \cite[Theorem A]{BFT2} that for each $k\geq 1$ the relation
\begin{eqnarray}
&&\int_{\Omega}|\nabla u|^pdx -\hk^p\int_{\Omega}\frac{|u|^p}{d^p}dx +
 \frac{1}{2}\hk^{p-1}\sum_{i=1}^{k-1}\int_{\Omega}\frac{|u|^p}{d^p}X_1^2\ldots X_i^2dx \nonumber \\
&\geq &c\int_{\Omega}\frac{|u|^p}{d^p}X_1^2\ldots X_k^{\gamma}dx
\label{eq:ggg}
\end{eqnarray}
is not valid for  $\gamma<2$; In addition, the best constant $c$ in (\ref{eq:ggg}) when $\gamma=2$ is equal to $\frac{1}{2}\hk^{p-1}$, for any $k=1,2, \ldots$.

A natural question is whether the constants appearing in (\ref{eq:t}) or  (\ref{eq:dos})
are  optimal.  Working towards this we consider the simplest  case (\ref{eq:l})
 (corresponding to  $p=2$, $k=0$).  Let  $\Omega=B$, be  the unit ball in $\R^N$,
and denote by $C_N$ the best constant of (\ref{eq:l}), that is
\beq
C_N =\inf_{u\in H^1_0(B)}\frac{\int_B|\nabla u|^2dx -\frac{1}{4}\int_B\frac{u^2}{d^2}dx}{\int_Bu^2dx}.
\label{eq:chr}
\eeq
We then show that in this case  the constant $k_2=\frac{N}{4}$ appearing in (\ref{eq:l}) is far from being optimal.
In particular we have:

{\bf Theorem B} {\em
For $N=3$, $C_3=\mu_2$, whereas for any $N \geq 2$ there holds:
 \beq 
C_N\geq \mu_2+\frac{(N-1)(N-3)}{4},
\label{1.20}
\eeq
where $ \mu_2\simeq 5.783$  is the best constant of inequality (\ref{eq:bv}).
}

It is remarkable that when $\Omega$ is a ball and  $N=3$  inequalities (\ref{eq:l}) and (\ref{eq:bv}) have the same
best constant.  For any  $N \geq 2$ the lower bound (\ref{1.20}) on $C_N$ improves
 the estimate $C_N \geq k_2 =\frac{N}{4}$.

To prove Theorem A we combine a  vector 
field approach  (cf [BFT]) along with ideas of  [HHL] or
[T]. It is worth noting that the ``mean distance'' method  of  Davies (cf [D1], [D2])
plays an essential role. For Theorem B after restricting to radial functions
we use suitable change of variables.

\setcounter{equation}{0}
\section{Preliminary inequalities}
In this section we will prove some auxiliary one-dimensional inequalities.
Throughout this section $b \leq \frac{\diam(\Omega)}{2}$ is a fixed positive constant.
 We have the following

\begin{lemma}
Let $\rho(t)=\min\{t,2b-t\}$. For any function $g\in C^1((0,b])$ there holds
\begin{eqnarray}
\ia && \int_0^{2b}|u'(t)|^pdt\geq \int_0^{2b}\{g'(\rho(t))-(p-1)|g(\rho(t))|^{\frac{p}{p-1}}
\}|u(t)|^pdt - 2 g(b) |u(b)|^p,  \nonumber \\
\ib &&
\int_0^{2b}|u'(t)|^pdt\geq \int_0^{2b}\{g'(\rho(t))-(p-1)|g(\rho(t))-g(b)|^{\frac{p}{p-1}}
\}|u(t)|^pdt
%\label{eq:rizi}
\end{eqnarray}
for all $u\in C^{\infty}_c(0,2b)$.
\label{lem:akr}
\end{lemma}
{\em Proof. } We first prove $\ia$. For $u\in C^{\infty}_c(0,2b)$ we have
\begin{eqnarray*}
\int_0^bg'(t)|u(t)|^pdt&=&g(b)|u(b)|^p-p\int_0^pg(t)|u|^pu'gdt\\
&\leq&g(b)|u(b)|^p+p\Bigl(\int_0^D|u'|^pdt\Bigr)^{\frac{1}{p}}
\Bigl(\int_0^b|g|^{\frac{p}{p-1}}|u|^pdt\Bigr)^{\frac{p-1}{p}}\\
&\leq&g(b)|u(b)|^p+\int_0^b|u'|^pdt+(p-1)\int_0^b|g|^{\frac{p}{p-1}}|u|^pdt,\\
\end{eqnarray*}
 hence
\[\int_0^b|u'(t)|^pdt\geq\int_0^b\{g'(t)-(p-1)|g(t)|^{\frac{p}{p-1}}\}
|u(t)|^pdt -g(b)|u(b)|^p.\]
A similar argument on $(b,2b)$ gives
\[\int_b^{2b}|u'(t)|^pdt\geq\int_b^{2b}\{g'(2b-t)-(p-1)|g(2b-t)|^{\frac{p}{p-1}}\}|u|^pdt
-g(b)|u(b)|^p\]
and  (i) follows by adding up the last two inequalities.

Part (ii) follows immediately  from (i) by using the function $g(x)-g(b)$  in the  place of $g(x)$.
$\hfill //$

In order to apply the above lemma we fix a positive integer $k$ and define the functions
\begin{eqnarray*}
&&\eta(t)=\sum_{i=1}^kX_1(t)\ldots X_i(t), \\
&& B(t)=\sum_{i=1}^kX_1^2(t)\ldots X_i^2(t),\quad t\in (0,1),
\end{eqnarray*}
where the $X_i$'s are given by (\ref{eq:flash}). It is easy to check that both $\eta$ and $B$ 
are increasing functions of $t$ with $\eta(0^+)=B(0^+)=0$ and $\eta(1^-)=B(1^-)=k$.
 We also  note that
\beq
\frac{1}{k}\eta^2(t) \leq B(t) \leq\eta^2(t) ,\quad t \in (0,1).
\label{eq:eb}
\eeq

 For $0<  b  \leq \frac{\diam(\Omega)}{2} \leq D$  we  define
 the following functions of $s \in (0,b)$:
\begin{eqnarray}
&&\hspace{-10pt}g(s)=-\Bigl(\frac{p-1}{p}\Bigr)^{p-1}s^{-(p-1)}\Bigl(1-\eta(s/D)-
a\eta^2(s/D)\Bigr)\label{eq:dal} \\
&&\hspace{-10pt}A(s)=g'(s)-(p-1)\Bigl|g(s)-g(b)\Bigr|^{\frac{p}{p-1}}
-\Bigl(\frac{p-1}{p}\Bigr)^ps^{-p}
-\frac{1}{2}\Bigl(\frac{p-1}{p}\Bigr)^{p-1}s^{-p}B(s/D).\nonumber
\end{eqnarray}
 Recall that $a$ is defined in (\ref{2.a}).
We then  have the following

\begin{lemma}
There exists $D_0 = D_0(k, p ,\diam(\Omega)) \geq  \frac{\diam(\Omega)}{2}$, such that 
for all $D \geq D_0$ there holds:
\begin{eqnarray}
\ia &&  1-\eta( \frac{\diam(\Omega)}{2D} )-a\eta^2(\frac{\diam(\Omega)}{2D}) \geq 0,  \nonumber  \\
\ib && g'(s)-\hk^ps^{-p}-\frac{1}{2}\hk^{p-1} s^{-p} B(s/D)
\geq (p-1)|g(s)|^{\frac{p}{p-1}}, \nonumber  \\
%\label{eq:dd} \\
\ic && \mbox{$A(s)$ is a decreasing function of $s\in (0,b)$.}\nonumber
\end{eqnarray}
\label{lem:dal}
For $p=2$, (ii) becomes equality. Also, for $p=2$, we can take as $D_0$ the unique solution
of $1=\eta( \frac{\diam(\Omega)}{2D_0} )$.
\end{lemma}
{\em Proof.} A straightforward calculation shows  that
\begin{equation}
\frac{d}{ds}\eta(s/D)=\frac{1}{s}\Bigl[\frac{B(s/D)}{2}+\frac{\eta^2(s/D)}{2}\Bigr].
\label{eq:k1}
\end{equation}
Setting $\Gamma(t)=tB'(t)$ we also have
\beq
\frac{d}{ds}B(s/D)=\frac{1}{s}\Gamma(s/D)>0;
\label{eq:k2}
\eeq
the positivity follows from the fact that $B(t)$ is an increasing function of $t$.

Since $\eta(t)$ is an increasing function of $t$ with $\eta(0)=0$, $\ia$ is immediate. 

We shall henceforth omit the argument $s/D$ from $\eta,B,\Gamma$ in the subsequent
formulas.
 We next  prove $\ib$. For $p=2$ an easy calculation shows that (ii) becomes equality.
For $p \neq 2$ the left hand side of $\ib$ is equal to
\begin{eqnarray}
&&g'(s)-\hk^ps^{-p}-\frac{1}{2}\hk^{p-1}s^{-p}B(s/D) =\hk^p(p-1)s^{-p}\times \nonumber \\
&&\hspace{0cm}\times \left[ 1-\frac{p\eta}{p-1}+(\frac{p}{2(p-1)^2}-\frac{ap}{p-1})\eta^2
 +\frac{ap}{(p-1)^2}\eta^3+\frac{ap}{(p-1)^2} \eta B \right]. 
\label{eq:meta}
\end{eqnarray}
On the other hand,
taking the Taylor expansion of $(1-t)^{\frac{p}{p-1}}$ about $t=0$, we see that
the right hand side of (ii) is written as (for $\eta$ small)
\begin{eqnarray}
&& \hk^p(p-1)s^{-p} (1-\eta-a\eta^2)^{\frac{p}{p-1}}  = 
\hk^p(p-1)s^{-p}\times  \label{eq:hop} \\
  && \times  \left[
1-\frac{p\eta}{p-1}-\frac{ap}{p-1}\eta^2 +\frac{p\eta^2}{2(p-1)^2}
+\frac{pa\eta^3}{(p-1)^2}
 +   \frac{p(p-2)\eta^3}{6(p-1)^3}+O(\eta^4) \right] .
 \nonumber
\end{eqnarray}
Comparing (\ref{eq:meta}) and (\ref{eq:hop}) we see that the corresponding right-hand sides agree
to order $O(\eta^2)$. Recalling (\ref{eq:eb}) and the choice of $a$ (cf (\ref{2.a})) we see that the cubic term in  (\ref{eq:meta}) is larger than the cubic term of (\ref{eq:hop}).  Hence (ii) is true 
 provided
$\eta$ is  small enough, which amounts to $D_0$ being large enough.

We now prove $\ic$. Note that $\ib$ implies that $g'$ is positive in $(0,b)$ if $D_0$ is large enough. Hence for $s \in (0,b)$ we have
\begin{eqnarray}
A'(s)&=&g''(s)+p[g(b)-g(s)]^{\frac{1}{p-1}}g'(s)+\hk^{p-1}(p-1)s^{-p-1}\nonumber\\
&&+\frac{1}{2}\hk^{p-1}ps^{-p-1}B-\frac{1}{2}\hk s^{-p-1}\Gamma\nonumber\\
&\leq&g''(s)+p|g(s)|^{\frac{1}{p-1}}g'(s)+\hk^{p-1}(p-1)s^{-p-1}\nonumber\\
&&+\frac{1}{2}\hk^{p-1}ps^{-p-1}B-\frac{1}{2}\hk s^{-p-1}\Gamma.
\label{eq:aris}
\end{eqnarray}
Using Taylor's expansion we have
\begin{eqnarray}
|g(s)|^{\frac{1}{p-1}}&=&\frac{p-1}{p}s^{-1}(1-\eta-a\eta^2)^{\frac{1}{p-1}}\nonumber\\
&=&\frac{p-1}{p}s^{-1}\Biggl\{1-\frac{1}{p-1}\eta -\Bigl[\frac{a}{p-1}+
\frac{p-2}{2(p-1)^2}\Bigr]\eta^2-\nonumber\\
&&-\Bigl[\frac{(p-2)a}{(p-1)^2}-\frac{(p-2)(3-2p)}{6(p-1)^3}\Bigr]\eta^3 +O(\eta^4)\Biggr\}.
\label{eq:snow}
\end{eqnarray}
From (\ref{eq:k1}), (\ref{eq:k2}), (\ref{eq:aris}) and (\ref{eq:snow}) we  obtain
\beq
A'(s)\leq(p-1)^2\hk^{p-1}s^{-p-1}
\Bigl\{\frac{p(p-2)}{6(p-1)^3}\eta^3-\frac{ap}{(p-1)^2}\eta B +O(\eta^4)\Bigr\}
\label{eq:adkr}
\eeq
From this and the fact  that
\[ \frac{1}{k}\eta^2 \leq B\leq\eta^2,\quad s\in (0,b)\]
we end up with
\beq
A'(s) \leq p \hk^{p-1}s^{-p-1} \eta^3 \Bigl\{ \frac{p-2}{6(p-1)} - \frac{a}{k} + O(\eta) \Bigr\}.
\label{222}
\eeq
To conclude the proof we distinguish various cases:\nl
(a) $1<p<2$. Then  $a=0$ and  it follows from (\ref{222})  that $A'(s)<0$ in $(0,b)$, provided
$D_0$ is chosen large enough.\nl
(b) $p=2$. Again $a=0$. A straightforward calculation shows that the right hand side of (\ref{eq:aris}) is identically equal to zero.
 The only restriction here comes from (i), whence the choice
of $D_0$. \nl
(c) $p>2$. Now  $a=\frac{(p-2) k}{3(p-1)}$ and the result follows again from  (\ref{222}). \nl
 This completes the proof.$\hfill //$

\setcounter{equation}{0}
\section{The Hardy inequality}

Throughout the rest of the paper we assume that $\Omega\subset\R^N$ is convex
and set $d(x)=\dist(x,\partial\Omega)$. 

Following \cite{HHL}, for $\omega\in S^{N-1}$ and $x\in\Omega$ we define
the following functions with values in $(0,+\infty]$:
\begin{eqnarray*}
\tau_{\omega}(x)&=&\inf\{s>0 \; |\; x+s\omega\not\in\Omega\} \\
\rho_{\omega}(x)&=&\min\{\tau_{\omega}(x),\tau_{-\omega}(x)\}\\
b_{\omega}(x)&=&\frac{1}{2}(\tau_{\omega}(x)+\tau_{-\omega}(x)).
\end{eqnarray*}

We denote by $dS(\omega)$ the standard measure on $S^{N-1}$ normalized so that
the total measure is one. Let $K_p>0$ be defined by
\beq
\int_{S^{N-1}}|v \cdot\omega|^pdS(\omega)=K_p  \, |v|^p ,\quad\forall v \in \R^N.
\label{eq:def}
\eeq
The constant $K_p$ is computable and with $k_p$ as in (\ref{1.kp}) we have
\begin{equation}
k_p=(p-1)\hk^pK_p^{-1}.
\label{eq:rep}
\end{equation}
We have the following
\blem
Assume that $\Omega$ is convex. Then for all $x  \in\Omega$ there holds
\beq
\int_{S^{N-1}}\rho_{\omega}^{-p}(x)dS(\omega)\geq K_p\,  d(x)^{-p}.
\eeq
\label{lem:mdf}
\elem
{\em Proof.} Let $y\in\partial\Omega$ be such that $|y-x|=d(x)$ and let
$P_y$ be the supporting hyper-plane through $y$ which is orthogonal to $y-x$.
We define the half-sphere
\[S^+=\{\omega\in S^{N-1}\; |\; \omega\cdot (y-x) >0\}\]
and for $\omega\in S^+$ define $\sigma_{\omega}(x)>0$ by requiring that
$x+\sigma_{\omega}(x)\omega\in P_y$, so that
\[\omega\cdot \frac{y-x}{|y-x|}=\frac{|y-x|}{\sigma_{\omega}(x)}.\]
The convexity of $\Omega$ implies that $\tau_{\omega}(x)\leq\sigma_{\omega}(x)$ and hence
\begin{eqnarray*}
\int_{S^{N-1}}\frac{1}{\rho_{\omega}(x)^p}dS(\omega)&\geq&
2\int_{S^+}\frac{1}{\tau_{\omega}(x)^p}dS(\omega) \\
&\geq&2\int_{S^+}\frac{1}{\sigma_{\omega}(x)^p}dS(\omega)\\
&=&\frac{2}{d(x)^{2p}}\int_{S^+}|(y-x)\cdot\omega|^pdS(\omega)\\
&=&\frac{K_p}{d(x)^p}\, ,
\end{eqnarray*}
as required.$\hfill //$

%%%%%%%%%   main theorem  %%%%%%%%%%%

We now give the proof of Theorem A.

{\em Proof of Theorem A.} Following \cite{HHL}
let us fix a direction $\omega\in S^{N-1}$ and let $\Omega_{\omega}$
be the orthogonal projection of $\Omega$ on the hyper-plane perpendicular to
$\omega$. For each $z\in\Omega_{\omega}$ we apply Lemma \ref{lem:akr} on
the segment defined by $z$ and $\omega$ and we then integrate over $z\in\Omega_{\omega}$.
We conclude that for any $u\in C^{\infty}_c(\Omega)$ there holds
\[\int_{\Omega}|\nabla u\cdot\omega|^pdx\geq\int_{\Omega}\left\{ g'(\rho_{\omega}(x))
-(p-1)\Bigl|g(\rho_{\omega}(x))-g(b_{\omega}(x))\Bigr|^{\frac{p}{p-1}}\right\}|u|^pdx,\]
Integrating over $\omega\in S^{N-1}$ and recalling definition (\ref{eq:def}) we obtain
\begin{eqnarray}
&&\int_{\Omega}|\nabla u|^pdx \geq K_p^{-1} \, \int_{\Omega}\int_{S^{N-1}}
\Bigl\{ g'(\rho_{\omega}(x))- \nonumber\\
&&\qquad\qquad -(p-1)\Bigl|g(\rho_{\omega}(x))-g(b_{\omega}(x))\Bigr|^{\frac{p}{p-1}}
\Bigr\}dS(\omega)|u|^pdx.
\label{eq:trip}
\end{eqnarray}
Now, let us choose $g$ as in (\ref{eq:dal}).
Since $\Omega$ is bounded Lemma \ref{lem:dal} implies the existence
of a  $D_0>0$ such that for $D \geq D_0$,  each of the functions
\begin{eqnarray*}
A_{\omega,x}(s)&:=&g'(s)-(p-1)\Bigl|g(s)-g(b_{\omega}(x))\Bigr|^{\frac{p}{p-1}}-\\
&&-\hk^ps^{-p}-\frac{1}{2}\hk^{p-1}s^{-p}B(s/D)
\end{eqnarray*}
-- defined for $s\in (0,b_{\omega}(x))$ -- is a decreasing function of
$s\in (0,b_{\omega}(x))$. In particular
$A_{\omega,x}(\rho_{\omega}(x))\geq A_{\omega,x}(b_{\omega}(x))$, i.e.
\begin{eqnarray*}
&&g'(\rho_{\omega}(x))-(p-1)\Bigl|g(\rho_{\omega}(x))-
g(b_{\omega}(x))\Bigr|^{\frac{p}{p-1}}\\
&\geq&\hk^p\rho_{\omega}(x)^{-p} +\frac{1}{2}\hk^{p-1}\rho_{\omega}(x)^{-p}
B(\rho_{\omega}(x)/D)+A_{\omega,x}(b_{\omega}(x)).
\end{eqnarray*}
Hence (\ref{eq:trip}) yields
\begin{eqnarray}
\int_{\Omega}|\nabla u|^pdx&\geq& K_p^{-1} \int_{\Omega}\int_{S^{N-1}}\Biggl\{
\hk^p\rho_{\omega}(x)^{-p} +\frac{1}{2}\hk^{p-1}\rho_{\omega}(x)^{-p}
B(\rho_{\omega}(x)/D)\nonumber\\
&&\qquad +g'(b_{\omega}(x))-\hk^pb_{\omega}(x)^{-p}\nonumber\\
&&\qquad -\frac{1}{2}\hk^{p-1}b_{\omega}(x)^{-p}B(b_{\omega}(x)/D)\Biggr\}dS(\omega)|u|^pdx.
\label{eq:sing}
\end{eqnarray}
We first estimate the first two terms of (\ref{eq:sing}).
For each $x\in\Omega$ and $\omega\in S^{N-1}$
there holds $B(\rho_{\omega}(x)/D)\geq B(d(x)/D)$,
and Lemma \ref{lem:mdf} yields
\begin{eqnarray}
&&K_p^{-1}  \, \int_{S^{N-1}}\Biggl\{\hk^p\rho_{\omega}(x)^{-p}
+\frac{1}{2}\hk^{p-1}\rho_{\omega}(x)^{-p}B(\rho_{\omega}(x)/D)\Biggr\}dS(\omega)\nonumber \\
&& \hspace{2cm}\geq\hk^p d(x)^{-p}+\frac{1}{2}\hk^{p-1}d(x)^{-p}B(d(x)/D),
\label{eq:adi}
\end{eqnarray}
for all $x\in\Omega$.
The remaining three terms in the right-hand side of (\ref{eq:sing}) are estimated using Lemma
\ref{lem:dal}(ii)
\begin{eqnarray*}
&&\hspace{-2cm}g'(b_{\omega}(x))-\hk^pb_{\omega}(x)^{-p}-\frac{1}{2}\hk^{p-1}b_{\omega}(x)^{-p}B( b_{\omega}(x)/D)
 \\
&\geq&(p-1)|g(b_{\omega}(x))|^{\frac{p}{p-1}} \nonumber \\
&=& \hk^p(p-1)(1-\eta_D-a\eta_D^2)^{\frac{p}{p-1}}  b_{\omega}(x)^{-p}.
\label{eq:stat}
\end{eqnarray*}
Combining this with  (\ref{eq:sing}) and  (\ref{eq:adi})  and recalling (\ref{eq:rep}) we obtain
\begin{eqnarray}
\int_{\Omega}|\nabla u|^pdx &\geq& \hk^p\int_{\Omega}\frac{|u|^p}{d^p}dx +
\frac{1}{2}\hk^{p-1}\int_{\Omega}\frac{|u|^p}{d^p}B(d/D)+\label{eq:tor}\\
&+  &k_p^{-1} \, 
(1-\eta_D-a\eta_D^2)^{\frac{p}{p-1}}   \int_{\Omega}\int_{S^{N-1}}\frac{1}{b_{\omega}(x)^p}
dS(\omega)|u|^pdx.\nonumber
\end{eqnarray}
We estimate the last integral using a variation of an argument of \cite{HHL}.
Elementary analysis shows that  $\min_{t>0}(1+t^N)/(1+t)^N =2^{-(N-1)}$ and therefore for $x\in\Omega$
\begin{eqnarray*}
2^{-\frac{(N-1)p}{N+p}}&\leq&\int_{S^{N-1}}\frac{(\tau_{\omega}(x)^N
+\tau_{-\omega}(x)^N)^{p/(N+p)}}
{(\tau_{\omega}(x)+\tau_{-\omega}(x))^{Np/(N+p)}}dS(\omega)\\
&\leq&\left(\int_{S^{N-1}}(\tau_{\omega}^N+\tau_{-\omega}^N)dS(\omega)
\right)^{\frac{p}{N+p}}
\left(\int_{S^{N-1}}\frac{1}{(\tau_{\omega}+\tau_{-\omega})^p}dS(\omega)
\right)^{\frac{N}{N+p}}\\
&=&2^{\frac{p-pN}{N+p}}
\left(\int_{S^{N-1}}\tau_{\omega}^NdS(\omega)\right)^{\frac{p}{N+p}}
\left(\int_{S^{N-1}}\frac{1}{b_{\omega}^p}dS(\omega)
\right)^{\frac{N}{N+p}},
\end{eqnarray*}
that is
\beq
\int_{S^{N-1}}\frac{1}{b_{\omega}(x)^p}dS(\omega)\geq
\left(\int_{S^{N-1}}\tau_{\omega}(x)^NdS(\omega)\right)^{-p/N}.
\label{eq:hot}
\eeq
The convexity of $\Omega$ implies $a_N\int_{S^{N-1}}\tau_{\omega}(x)^NdS(\omega)=
|\Omega|$. Hence the proof
is concluded by combining (\ref{eq:tor}) and (\ref{eq:hot}).$\hfill //$

{\bf Remark} We note that inequality (\ref{eq:trip}) can be used to obtain
Hardy type  inequalities for non convex domains as in [HHL], [T].

\setcounter{equation}{0}
\section{On the best constant for $p=2$}

In this section we will prove Theorem B. We recall that 
 $C_N$ is  the best constant of inequality (\ref{eq:l}), in case $\Omega$  is  a  ball, 
 defined by:
\beq
C_N =\inf_{u\in H^1_0(B)}\frac{\int_B|\nabla u|^2dx -\frac{1}{4}\int_B\frac{u^2}{d^2}dx}{\int_Bu^2dx}.
\label{4.1}
\eeq
We first establish
\blem
The infimum in (\ref{4.1})  remains the same if it is taken over all radially symmetric
functions $u=u(r)\in H^1_0(B)$.
\label{lem.41}
\elem
{\em Proof:}
We may assume that $\Omega$ is the unit ball.
 Let us denote by $\tilde{C}_N$ the infimum over radial
functions. Clearly $\tilde{C}_N\geq C_N$.
Suppose now that 
%$u\in C^{\infty}_c(B)$    
$ u \in  H^1_0(B)$ and let
\[
u(x)=u_0(r)+\sum_{m=1}^{\infty}f_m(\sigma)u_m(r), \quad \quad r=|x|,
\]
be its decomposition into spherical harmonics; here $u_m$ are radially symmetric functions in
$H^1_0(B)$   
% $C^{\infty}_c(B)$
   and $f_m$ are orthonormal in $L^2(S^{N-1})$  eigenfunctions
of the Laplace-Beltrami operator on $\{|x|=1\}$, with corresponding
eigenvalues $c_m=m(N-2+m)$, $m\geq 1$. It is easily seen that
\begin{equation}
\int_{B}|\nabla u|^2dx =\int_{B}(|\nabla u_0|^2 dx+\sum_{m=1}^{\infty}\int_{B}(|\nabla u_m|^2
+\frac{c_m}{|x|^2}u_m^2)dx,
\label{eq:oma}
\end{equation}
and hence
\begin{eqnarray*}
\int_B(|\nabla u|^2-\frac{u^2}{d^2})dx &=&\int_B \{|\nabla u_0|^2 -\frac{u_0^2}{4(1-|x|)^2}\}dx +\\
&&\qquad +\sum_{m=1}^{\infty}\int_B\{(|\nabla u_m|^2  +(\frac{c_m}{|x|^2}-\frac{1}{4(1-|x|)^2}) u_m^2\}dx\\
&\geq&\tilde{C}_N\int_B u_0^2 +\tilde{C}_N \sum_{m=1}^{\infty}\int_B u_m^2dx \\
&=&\tilde{C}_N \int_Bu^2dx.
\end{eqnarray*}
This implies $C_N\geq\tilde{C}_N$ and the Lemma  is proved. $\hfill //$

{\em Proof of Theorem B:} By the previous Lemma
 we restrict attention to radially symmetric functions.
 Let $u = u(r)\in C^{\infty}_c(B)$ be a  radial function  and define $v$ by
\[ u(r)=r^{-\frac{N-1}{2}}(1-r)^{1/2}v(r),\qquad r\in (0,1).\]
Then $v(0)=v(1)=0$. We compute
\begin{eqnarray*}
\frac{1}{N a_N}
\int_B|\nabla u|^2dx&=&\int_0^1 (u')^2r^{N-1}dr \\
&=&\int_0^1(1-r)\Bigl(-\frac{(N-1)v}{2r}-\frac{v}{2(1-r)}+v'\Bigr)^2dr
\end{eqnarray*}
Using integration by parts for the terms involving $vv'=(v^2)'/2$ we conclude after some
 simple calculations that
\begin{eqnarray*}
\frac{1}{N a_N} \left(
\int_B|\nabla u|^2dx -\frac{1}{4}\int_B\frac{u^2}{d^2}dx \right) &=&\int_0^1 (1-r)(v')^2dr +\frac{(N-1)(N-3)}{4}\int_0^1\frac{1-r}{r^2}v^2dr\\
&\geq&\int_0^1 (1-r)(v')^2dr +\frac{(N-1)(N-3)}{4}\int_0^1(1-r)v^2dr.
\end{eqnarray*}
But (cf [BV, Section 4])
\[\inf_{v(0)=v(1)=0}\frac{\int_0^1(1-r)(v')^2dr}{\int_0^1(1-r)v^2dr}=\inf_{v(0)=v(1)=0}\frac{\int_0^1r(v')^2dr}{\int_0^1rv^2dr}=\mu_2,
\]
and estimate (\ref{1.20}) of Theorem B follows.

To prove that $C_3=\mu_2$, let us define
\[ u_{\epsilon}(r)=r^{-1}(1-r)^{\frac{1}{2}+\epsilon}w(1-r),\qquad r\in (0,1),\]
where $\epsilon>0$ and $w(|x|)$ is the first eigenfunction of the Dirichlet Laplacian for the unit
disk in $\R^2$. Then 
\[u_{\epsilon}'(r)=-r^{-1}(1-r)^{\frac{1}{2}+\epsilon}\Bigl\{\frac{w}{r}
+(\frac{1}{2}+\epsilon)\frac{w}{1-r} +w'\Bigr\}\]
and hence $u_{\epsilon}\in H^1_0(B)$ and
\begin{eqnarray*}
&&\int_0^1(u_{\epsilon}')^2r^{N-1}dr=\int_0^1(1-r)^{1+2\epsilon}\left\{\frac{w^2}{r^2}
 +(\frac{1}{2}+\epsilon)^2\frac{w^2}{(1-r)^2}
+(w')^2 \right. \\
&& \left. +(1+2\epsilon)\frac{w^2}{r(1-r)} +\frac{2ww'}{r}
+\frac{(1+2\epsilon)ww'}{1-r}
\right\}dr \qquad \mbox{ (where $w=w(1-r)$)}
\end{eqnarray*}
To handle the terms containing $ww'$ we
integrate by parts: the boundary terms are equal to zero and making the change
of variables $s=1-r$ we eventually obtain
\[\int_0^1\Bigl((u_{\epsilon}'(r))^2-\frac{1}{4}\frac{u_{\epsilon}^2(r)}{(1-r)^2}\Bigr)r^2dr
=\int_0^1 s^{1+2\epsilon}\Bigl( (w'(s))^2-\epsilon^2\frac{w^2(s)}{s^2}\Bigr)ds.\]
Now, there holds
\[ \epsilon^2\int_0^1  s^{-1+2\epsilon}w^2 ds\longrightarrow 0,\qquad
\mbox{ as }\epsilon\to 0,\]
hence
\begin{eqnarray*}
\lim_{\epsilon\to 0}\frac{\int_B\Bigl(|\nabla u_{\epsilon}|^2-\frac{u_{\epsilon}^2}{
4d^2}\Bigr)dx}{\int_Bu_{\epsilon}^2dx}&=&\lim_{\epsilon\to 0}
\frac{\int_0^1 (w')^2 s^{1+2\epsilon}ds}{\int_0^1 w^2s^{1+2\epsilon}ds}\\
&=&\frac{\int_0^1 (w')^2 s\,ds}{\int_0^1 w^2s\,ds}\\
&=&\mu_2.
\end{eqnarray*}
It follows that $\tilde{C}_3\leq \mu_2$; in view of (\ref{1.20})
and Lemma \ref{lem.41} we conclude that  $C_3=\mu_2$.

%%%%%%%%%%%%%%%%%%%%%%%%%%%%%%%%%%%

%%%%%%%%%%%%   bibliography  %%%%%%%%%%%%

%%%%%%%%%%%%%%%%%%%%%%%%%%%%%%%%%%%

\end{document}